# MOD 2 REPRESENTATIONS OF ELLIPTIC CURVES

K. RUBIN AND A. SILVERBERG

ABSTRACT. Explicit equations are given for the elliptic curves (in characteristic $\neq 2, 3$) with mod 2 representation isomorphic to that of a given one.

## 1. INTRODUCTION

If $N$ is a positive integer and $E$ is an elliptic curve defined over a field $F$, one can ask for a description of the set of elliptic curves whose mod $N$ representation (of the absolute Galois group) is symplectically isomorphic to that of $E$ (see [2]). For $N = 3$, 4, and 5, we gave explicit equations in [3] and [5]. The case $N = 1$ is trivial, and when $N \geq 7$ the set in question is always finite and the situation is quite different from the ones we consider. In [4] we gave a description for $N = 6$ (but did not give explicit equations).

This note, which can be viewed as a footnote to those papers, deals with the easier case $N = 2$. Note that since there is only one nondegenerate alternating pairing on $\mathbf{Z}/2\mathbf{Z} \times \mathbf{Z}/2\mathbf{Z}$, isomorphic and symplectically isomorphic are the same for mod 2 representations. Theorem 1 gives explicit equations for the family of elliptic curves whose mod 2 representation is isomorphic to that of a given one. Given two elliptic curves, Corollary 2 gives an easy way to determine whether or not their mod 2 representations are isomorphic. The proofs are given in §2. In §3 we give a different approach, using the algorithm from [3].

If $F$ is a field, let $F^{\text{sep}}$ denote a separable closure of $F$ and let $G_F = \text{Gal}(F^{\text{sep}}/F)$. If $E$ is an elliptic curve over $F$, let $j(E)$ denote its $j$-invariant, let $\Delta(E)$ denote its discriminant, and let $E[2]$ denote the $G_F$-module of 2-torsion points on $E$.

**Theorem 1.** *Suppose $F$ is a field of characteristic different from 2 and 3, and $E : y^2 = x^3 + ax + b$ is an elliptic curve over $F$. If $u, v \in F$, let $\mathcal{E}_{u,v}$ denote the curve*

$$y^2 = x^3 + 3(3av^2 + 9buv - a^2u^2)x + 27bv^3 - 18a^2uv^2 - 27abu^2v - (2a^3 + 27b^2)u^3.$$

*If $E'$ is an elliptic curve over $F$, and $E'[2] \cong E[2]$ as $G_F$-modules, then $E'$ is isomorphic to $\mathcal{E}_{u,v}$ for some $u, v \in F$. Conversely, if $u, v \in F$ and $\mathcal{E}_{u,v}$ is nonsingular, then $\mathcal{E}_{u,v}[2] \cong E[2]$ as $G_F$-modules,*

$$j(\mathcal{E}_{u,v}) = \frac{(3av^2 + 9buv - a^2u^2)^3 j(E)}{27a^3(v^3 + au^2v + bu^3)^2}, \quad \text{and} \quad \Delta(\mathcal{E}_{u,v}) = 3^6(v^3 + au^2v + bu^3)^2 \Delta(E).$$

**Corollary 2.** *Suppose $F$ is a field of characteristic different from 2 and 3, and $E : y^2 = x^3 + ax + b$ is an elliptic curve over $F$. Let*

$$C(u,v) = \frac{(3av^2 + 9buv - a^2u^2)^3}{27a^3(v^3 + au^2v + bu^3)^2}.$$

*Suppose $E'$ is an elliptic curve over $F$. If $j(E') \neq 0, 1728$, and for some $(u,v) \in \mathbf{P}^1(F)$ we have*

(i) $\quad\dfrac{j(E')}{j(E)} = C(u,v) \quad$ if $a \neq 0$, $\quad$ or

(ii) $\quad\dfrac{j(E')}{j(E) - 1728} = \dfrac{-4C(u,v)a^3}{27b^2} \quad$ if $b \neq 0$,

then $E'[2] \cong E[2]$. Conversely, if $E'[2] \cong E[2]$, then there is a point $(u,v) \in \mathbf{P}^1(F)$ such that $j(E')$ satisfies (i) if $a \neq 0$ and (ii) if $b \neq 0$.

We thank the NSF and NSA for financial support. Silverberg thanks AIM and the UC Berkeley math department for their hospitality.

## 2. Proofs

**Lemma 3.** *Suppose $F$ is a field and $\varphi(x) \in F[x]$ is a polynomial with no multiple roots. Let $\Psi_\varphi$ denote the set of roots of $\varphi$.*
  (i) *There is a $G_F$-equivariant bijection $\Psi_\varphi \xrightarrow{\sim} \mathrm{Hom}_{F-\mathrm{algebra}}(F[x]/(f(x)), F^{\mathrm{sep}})$.*
  (ii) *The $F$-algebra of $G_F$-equivariant maps from $\Psi_\varphi$ to $F^{\mathrm{sep}}$ is isomorphic to $F[x]/(f(x))$.*

*Proof.* Assertion (i) is clear, and (ii) follows from Lemma 5 on p. A.V.75 of [1]. □

**Lemma 4.** *Suppose $E : y^2 = f(x)$ and $E' : y^2 = g(x)$ are elliptic curves over a field $F$ with $f(x), g(x) \in F[x]$ of degree 3. Then $E[2] \cong E'[2]$ as $G_F$-modules if and only if $F[x]/(f(x)) \cong F[x]/(g(x))$ as $F$-algebras.*

*Proof.* We apply Lemma 3 with $\varphi = f$ and $g$. Since the roots of $f$ are the $x$-coordinates of the elements of $E[2] - 0$, there is a $G_F$-equivariant bijection $\Psi_f \xrightarrow{\sim} E[2] - 0$. Similarly we have $\Psi_g \xrightarrow{\sim} E'[2] - 0$. Thus by Lemma 3, $F[x]/(f(x)) \cong F[x]/(g(x))$ as $F$-algebras if and only if $E[2] - 0 \cong E'[2] - 0$ as $G_F$-sets. Since every bijection $E[2] - 0 \xrightarrow{\sim} E'[2] - 0$ extends to a group isomorphism $E[2] \xrightarrow{\sim} E'[2]$, the lemma follows. □

*Proof of Theorem 1.* Write $f(x) = x^3 + ax + b$, so $E$ is the elliptic curve $y^2 = f(x)$, and let $E'$ be an elliptic curve $y^2 = g(x) = x^3 + \alpha x + \beta$ with $\alpha, \beta \in F$.

Suppose $E[2] \cong E'[2]$ as $G_F$-modules. By Lemma 4, there is an isomorphism of $F$-algebras $\phi : F[z]/(g(z)) \xrightarrow{\sim} F[x]/(f(x))$. Let $\phi(z) = 3ux^2 + 3vx + w$ with $u, v, w \in F$. (The extra factors of $3$ remove denominators which would otherwise occur in the equation for $\mathcal{E}_{u,v}$ and the formulas below.) The trace of $\phi(z)$ acting by multiplication on $F[x]/(f(x))$ is $3w - 6au$, but the trace of $z$ acting by multiplication on $F[z]/(g(z))$ is zero. Since $\phi$ is an isomorphism, we must have $w = 2au$. Using this, the characteristic polynomial of $\phi(z)$ acting on $F[x]/(f(x))$ is

$$h(T) = T^3 + 3(3av^2 + 9buv - a^2u^2)T + 27bv^3 - 18a^2uv^2 - 27abu^2v - (2a^3 + 27b^2)u^3.$$

Again, since $\phi$ is an isomorphism, we conclude that $h(T) = g(T)$, i.e., $E'$ is $\mathcal{E}_{u,v}$ as desired.

Conversely, suppose that $u, v \in F$ are such that

$$\alpha = 3(3av^2 + 9buv - a^2u^2), \quad \beta = 27bv^3 - 18a^2uv^2 - 27abu^2v - (2a^3 + 27b^2)u^3.$$

Then working backwards through the argument above, one can show that the map $z \mapsto 3ux^2 + 3vx + 2au$ induces a homomorphism $\phi : F[z]/(g(z)) \to F[x]/(f(x))$. The determinant of $\phi$ with respect to the bases $\{1, z, z^2\}$ and $\{1, x, x^2\}$ is $27(v^3 + au^2v + bu^3)$. However, the discriminant of $g$ is $3^6(4a^3 + 27b^2)(v^3 + au^2v + bu^3)^2$. Since $E'$ is an elliptic curve, the discriminant of $g$ must be nonzero, and hence the determinant of $\phi$ is nonzero so $\phi$ is an isomorphism. By Lemma 4, it follows that $E[2] \cong E'[2]$ as $G_F$-modules.

The formulas for the $j$-invariant and the discriminant are immediate. $\square$

*Proof of Corollary 2.* If $u, v \in F$ are such that $j(E')$ satisfies (i) or (ii), then $\mathcal{E}_{u,v}$ is nonsingular (by the computation of its discriminant in Theorem 1) and $j(E') = j(\mathcal{E}_{u,v})$. If $j(E') \neq 0, 1728$, then $E'$ is a quadratic twist of $\mathcal{E}_{u,v}$. Therefore using Theorem 1, we have $E'[2] \cong \mathcal{E}_{u,v}[2] \cong E[2]$. Conversely, if $E'[2] \cong E[2]$, then by Theorem 1 we can find $u, v \in F$ such that $E' \cong \mathcal{E}_{u,v}$. By Theorem 1 we have (i) and (ii). $\square$

## 3. A different method

Applying the method of [3] (see also §3 of [5]) to the case $N = 2$, one again obtains explicit equations for the family of elliptic curves with mod 2 representation isomorphic to that of $E$. We show below how the algorithm works in this case. Suppose $F$ is a field with $\operatorname{char}(F) \neq 2, 3$, and $E : y^2 = x^3 + ax + b$ is an elliptic curve over $F$. Note that mod 2 representations do not change under quadratic twist. Every elliptic curve $E'$ over $F$ such that the $G_F$-action on $E'[2]$ is trivial is a quadratic twist of
$$A_\lambda : y^2 = x(x-1)(x-\lambda)$$
with $\lambda \in F - \{0, 1\}$. Putting $A_\lambda$ in Weierstrass form we obtain
$$E_\lambda : y^2 = x^3 + a_4(\lambda)x + a_6(\lambda),$$
where
$$a_4(\lambda) = -\frac{1}{3}(\lambda^2 - \lambda + 1), \quad a_6(\lambda) = -\frac{1}{27}(2\lambda^3 - 3\lambda^2 - 3\lambda + 2).$$
The algorithm in §3 of [3] shows that the equations we are looking for are of the form
$$(1) \qquad dy^2 = x^3 + a(t)x + b(t)$$
with
$$d \in F, \quad a(t) = \mu^{-2}(\gamma t + 1)^2 a_4(A(t)), \quad \text{and} \quad b(t) = \mu^{-3}(\gamma t + 1)^3 a_6(A(t)),$$
where $u_0$ satisfies $j(E_{u_0}) = j(E)$, $\mu$ satisfies
$$a_4(u_0) = a\mu^2 \quad \text{and} \quad a_6(u_0) = b\mu^3,$$
and
$$A(t) = \frac{\alpha t + u_0}{\gamma t + 1}$$
with $\alpha$ and $\gamma$ chosen so that $a(t), b(t) \in F[t]$.

If $ab \neq 0$, let $j = j(E)$ and let $u_0$ be a root of the numerator (as a polynomial in $\lambda$) of
$$j(E_\lambda) - j = \frac{256 - 768\lambda + (1536 - j)\lambda^2 + (2j - 1792)\lambda^3 + (1536 - j)\lambda^4 - 768\lambda^5 + 256\lambda^6}{\lambda^2(\lambda-1)^2}.$$

Let
$$\mu = \frac{a_6(u_0)a}{a_4(u_0)b} = \frac{(2u_0^3 - 3u_0^2 - 3u_0 + 2)a}{9(u_0^2 - u_0 + 1)b} \in (F^{\text{sep}})^\times,$$

$$\alpha = \frac{3(u_0 - 2)\mu^3 b}{u_0(u_0 - 1)}, \qquad \gamma = \frac{3(2u_0 - 1)\mu^3 b}{u_0(u_0 - 1)} \quad \in F^{\text{sep}}.$$

With these values, equation (1) becomes

$$dy^2 = x^3 + a(1 + (J-1)t^2)x + b(1 + 3t - 3(J-1)t^2 - (J-1)t^3),$$

where

$$J = \frac{j(E)}{1728} = \frac{4a^3}{4a^3 + 27b^2}.$$

For $d \in F$ and $t \in \mathbf{P}^1(F)$, this gives the elliptic curves over $F$ with mod 2 representation isomorphic to that of $E$, when $ab \neq 0$.

Similarly, if $b = 0$, then

$$j(E_\lambda) - j(E) = \frac{64\left(-2 + \lambda\right)^2 (1 + \lambda)^2 (-1 + 2\lambda)^2}{(-1 + \lambda)^2 \lambda^2}.$$

With $u_0 = 2$, $\mu = 1/\sqrt{-a}$, $\alpha = 0$, and $\gamma = 3\sqrt{-a}$, equation (1) becomes

$$dy^2 = x^3 + a(1 - 3at^2)x + 2a^2 t(1 + at^2).$$

If $a = 0$, then

$$u_0 = \frac{1 + \sqrt{-3}}{2}, \quad \mu = \frac{-1}{b^{1/3}\sqrt{-3}}, \quad \alpha = \frac{b^{1/3}(1 - \sqrt{-3})}{2}, \quad \text{and} \quad \gamma = b^{1/3}$$

yield the equation

$$dy^2 = x^3 + 3btx + b(1 - bt^3).$$

## References

[1] N. Bourbaki, Algebra II, Springer, Berlin, 1990.
[2] B. Mazur, *Rational isogenies of prime degree*, Invent. Math. **44** (1978), 129–162.
[3] K. Rubin, A. Silverberg, *Families of elliptic curves with constant mod p representations*, in Conference on Elliptic Curves and Modular Forms, Hong Kong, December 18–21, 1993, Intl. Press, Cambridge, Massachusetts, 1995, pp. 148–161.
[4] ______, *Mod 6 representations of elliptic curves*, To appear in Proc. Symp. Pure Math., AMS, Providence.
[5] A. Silverberg, *Explicit families of elliptic curves with prescribed mod N representations*, in Modular Forms and Fermat's Last Theorem, eds. Gary Cornell, Joseph H. Silverman, Glenn Stevens, Springer, Berlin (1997), 447–461.

DEPARTMENT OF MATHEMATICS, STANFORD UNIVERSITY, STANFORD, CA,
DEPARTMENT OF MATHEMATICS, OHIO STATE UNIVERSITY, COLUMBUS, OHIO 43210
*E-mail address*: rubin@math.stanford.edu

DEPARTMENT OF MATHEMATICS, OHIO STATE UNIVERSITY, 231 W. 18 AVENUE, COLUMBUS, OHIO 43210
*E-mail address*: silver@math.ohio-state.edu